\documentclass[12pt]{article}
\usepackage{latexsym}
\usepackage{amssymb}
\usepackage{graphicx}

\newtheorem{Theorem}{Theorem}

\newtheorem{Proposition}{Proposition}

\newtheorem{Lemma}{Lemma}
\newtheorem{Corollary}{Corollary}

\newtheorem{Example}{Example}
\newtheorem{Algorithm}{Algorithm}

\def \ep{\hbox{ }\hfill$\Box$}

\def \ra{\rightarrow}
\def\reff#1{{\rm(\ref{#1})}}

\begin{document}
\title{Convergence of a Second Order Markov Chain}

\author{
Shenglong Hu \thanks{Corresponding author. Email: Tim.Hu@connect.polyu.hk. Department of
Applied Mathematics, The Hong Kong Polytechnic University, Hung Hom,
Kowloon, Hong Kong. Tel: 852-2766 4368, Fax: 852-2362 9045.},\hspace{4mm} Liqun Qi \thanks{Email:
maqilq@polyu.edu.hk. Department of Applied Mathematics, The Hong
Kong Polytechnic University, Hung Hom, Kowloon, Hong Kong. This
author's work was supported by the Hong Kong Research Grant
Council. (Grant No. PolyU 501808, 501909, 502510 and 502111).} }

\date{\today} \maketitle

\begin{abstract}
In this paper, we consider convergence properties of a second order
Markov chain. Similar to a column stochastic matrix is associated to a
Markov chain, a so called {\em transition probability tensor} $P$ of
order $3$ and dimension $n$ is associated to a second order Markov
chain with $n$ states. For this $P$, define $F_P$ as
$F_P(x):=Px^{2}$ on the $n-1$ dimensional standard simplex
$\Delta_n$. If $1$ is not an eigenvalue of $\nabla F_P$ on
$\Delta_n$ and $P$ is irreducible, then there exists a unique fixed
point of $F_P$ on $\Delta_n$. In particular, if every entry of $P$
is greater than $\frac{1}{2n}$, then $1$ is not an eigenvalue of
$\nabla F_P$ on $\Delta_n$. Under the latter condition, we further
show that the second order power method for finding the unique fixed
point of $F_P$ on $\Delta_n$ is globally linearly convergent and the
corresponding second order Markov process is globally $R$-linearly
convergent.   \vspace{3mm}

\noindent {\bf Key words:}\hspace{2mm} Nonnegative tensor, transition probability tensor, second order Markov Chain
\vspace{3mm}

\noindent {\bf AMS subjection classifications (2010):}\hspace{2mm}15A18; 15A69; 60J10; 60J22; 65F15
\end{abstract}

\section{Introduction}
\hspace{4mm} Markov chain serves as a fundamental tool for diverse
applications \cite{lr,r,s}. As a generalization, higher order Markov
chain can be used to describe many phenomena in science and engineering, e.g., bioinformatics, genome, speech/text recogonition, please refer to \cite{cn,ln} and references therein.
Compared with the sophisticated development of Markov chain based on
stochastic matrices, research on higher order Markov chain based on
transition probability tensors, is just on the way
\cite{cn,ln,lny,nqz}. Nevertheless, the recent progress in numerical multilinear algebra, especially in tensor decomposition \cite{kn,kb} and spectral theory of tensors \cite{l,q}, introduces many new tools for this topic.

An $m$-th order $n$ dimensional Markov chain is basically
characterized by its associated nonnegative tensor $P$ which is an
$(m+1)$-th order $n$ dimensional tensor with entries $p_{ij_1\ldots
j_{m}}$ for all $i,j_1,\ldots,j_{m}=1,\ldots,n$ satisfying:
\begin{eqnarray*}
0\leq p_{ij_1\cdots j_{m}}=\mbox{Prob}(X_t=i\;|\;X_{t-1}=j_1,\ldots,X_{t-m}=j_{m})\leq 1.
\end{eqnarray*}
Here $\{X_t,t=0,1,\ldots\}$ represents the stochastic process that takes on $n$ states $\{1,\ldots,n\}$.
Obviously,
\begin{eqnarray}\label{int-1}
\sum_{i=1}^np_{ij_1\cdots j_{m}}=1
\end{eqnarray}
for any $j_1,\ldots,j_{m}=1,\ldots,n$.

An $(m+1)$-th order $n$ dimensional nonnegative tensor $P$ that
satisfies \reff{int-1} is called a {\it transition probability
tensor}. In this paper, we consider convergence properties of a
second order Markov chain which is associated to a third order $n$
dimensional transition probability tensor $P$. Given such a tensor
$P$, we specially consider the sequence of state distribution
vectors generated by the second order Markov process: with two
initial state distribution vectors $x^{(0)},
x^{(1)}\in\Delta_n:=\{x\in\Re^n_+\;|\;\sum_{i=1}^nx_i=1\}$, the
sequence is generated as
\begin{eqnarray}\label{int-2}
x^{(s)}:=Px^{(s-1)}x^{(s-2)},\;\forall s=2,3,\ldots,
\end{eqnarray}
where $Px^{(s-1)}x^{(s-2)}$ is an $n$-vector whose $i$-th element is
\begin{eqnarray*}
\sum_{j,k=1}^np_{ijk}x^{(s-1)}_jx^{(s-2)}_k
\end{eqnarray*}
for all $i\in\{1,\ldots,n\}$.
If the sequence $\{x^{(k)}\}$ converges to $x^*$, then we call $x^*$ the stationary probability distribution of the second order Markov chain. Obviously, in this situation,
\begin{eqnarray}\label{int-3}
x^*=Px^*x^*=:P(x^*)^2
\end{eqnarray}
which is closely related to Z-eigenvalues of tensors introduced in \cite{q}.

For the convenience of the subsequent analysis, define a nonlinear map $F_P:\Re^n\ra\Re^n$ associated to $P$ as:
\begin{eqnarray}\label{eq-1}
(F_P(x))_i=\sum_{j,k=1}^np_{ijk}x_{j}x_{k}
\end{eqnarray}
for all $i\in\{1,\ldots,n\}$ and $x\in\Re^n$, and denote by $\nabla F_P(x)$ the Jacobian matrix of the map $F_P$ at $x$. Obviously, $\sum_{i=1}^n(F_P(x))_i=1$ for all $x\in\Delta_n$.
Essentially, stationary probability distribution of the second order Markov chain associated to $P$ in \reff{int-3}, whenever it exists, is a fixed point of $F_P$ on $\Delta_n$.

Very recently, under mild conditions, some results like the
uniqueness of $x^*$ in \reff{int-3} and the linear convergence of
the power method for finding such a unique $x^*$ were established in
\cite{ln,lny}. Unlike its counterpart
of $m=1$ \cite{s}, the convergence of the original Markov process
\reff{int-2} could not be deduced directly from the convergence of
the power method for finding $x^*$ in \reff{int-3}. Hence, in this
paper, we mainly consider this problem and give an answer.

We give a uniqueness property of the fixed
point of $F_P$ on $\Delta_n$ by using fixed point index theory in
Section 2. For the sake of the clarification, the detailed proof is given in Appendix at the end of this paper. We establish globally linear convergence of the second
order power method for finding the unique fixed point in Section 3.
Globally $R$-linear convergence of the second order Markov process
\reff{int-2} is proved in Section 4. Some intuitive numerical examples
are given in the last section.

\section{The uniqueness property}
\hspace{4mm} In this section, we discuss the uniqueness of the fixed point of $F_P$ on $\Delta_n$ for a given transition probability tensor $P$, which is parallel to the strong Perron-Frobenius theorem for primitive stochastic matrices \cite{s}. In the first subsection, we give a general result, and we discuss more on third order transition probability tensors in the second subsection.

\subsection{General case}
\hspace{4mm} We give the general result for $m$-th order $n$ dimensional transition probability tensors. The following concept is useful: an $m$-th order $n$ dimensional
nonnegative tensor $P$ is called {\em reducible} if there exists a
nonempty proper index subset $I\subset \{1,\ldots,n\}$ such that
\begin{eqnarray}\label{red}
P_{ij_1\ldots j_{m-1}}=0,\quad \forall i\in I,\quad \forall
j_1,\ldots,j_{m-1}\notin I.
\end{eqnarray}
If $P$ is not reducible, then $P$ is called {\em irreducible} \cite{cpz}. Obviously, $P$ is irreducible if it is a positive tensor. The concepts of
relative interior and relative boundary of a set are in the usual
sense \cite{rok}.

In the subsequent anslysis, $P$ is assumed
to be an $m$-th order $n$ dimensional transition probability tensor.
In order to accomplish the proof, we first introduce briefly the
concept of fixed point index, see \cite{gl} for a comprehensive
discussion. The theory established in \cite{gl} applies to general
Banach spaces, while we present it in the finitely dimensional cases
to match our problem.  Intuitive speaking, the theory of fixed point
index is just a generalization of the degree theory: for the degree
theory discusses the fixed points of a map in an open set; while,
the theory of fixed point index discusses the fixed points of a map
in a relative open set. As our domain $\Delta_n$ here is relative
open in $\Re^n$, we need the theory of fixed point index to give
rigorous proof. Here is the analogue of fixed point index to degree
theory.

\begin{Theorem}\label{thm-2}
Suppose that $U\in\Re^n$ is a bounded relatively open subset with
closure $\overline{ U}$ and $f:\overline{ U}\ra\Re^n$ is continuous.
If $f$ has no fixed points on the relative boundary $\partial U$ of
$U$, then there exists an integer function, denoted as $Ind(f,U)$,
satisfying the following properties:
\begin{itemize}
\item [(A1)] Normality: $Ind(f,U)=1$ if $f(x)\equiv x_0\in U$ for all $x\in \overline{ U}$.
\item [(A2)] Additivity: If $U_1$ and $U_2$ are disjoint open subsets relative to $U$ and $f$ has no fixed points
on $\overline{ U}\setminus (U_1\cup U_2)$, then
$Ind(f,U)=Ind(f,U_1)+Ind(f,U_2)$.
\item [(A3)] Homotopy Invariance: If $H:[0,1]\times \overline{ U}\ra\Re^n$ is continuous and $H(t,\cdot)$
has no fixed points on $\partial U$ for any $t\in[0,1]$, then
$Ind(H(t,\cdot),U)$ is constant for all $t\in[0,1]$.
\item [(A4)] Excision: If $U_0$ is an open subset relative to $U$ and $f$ has no fixed points in
$\overline{ U}\setminus U_0$, then $Ind(f,U)=Ind(f,U_0)$.
\item [(A5)] Solution: If $Ind(f,U)\neq 0$, then $f$ has at least one fixed point in $U$.
\end{itemize}
In addition, $Ind(f,U)$ is uniquely defined.
\end{Theorem}

Now, we give the main result in this section.
\begin{Theorem}\label{thm-3}
Suppose that $P$ is an $m$-th order $n$ dimensional irreducible transition probability tensor and $1$ is not an eigenvalue of $\nabla F_P$ for all $x\in\mbox{rel}(\Delta_n)$, the relative interior of $\Delta_n$. Then, there is a unique $x\in\mbox{rel}(\Delta_n)$ such that $F_P(x)=x$.
\end{Theorem}

\noindent {\bf Proof.} The proof consists of four parts as follows.

{\bf I}. $F_P$ has no fixed points on $\partial\Delta_n$.

Suppose not, then there exists $x\in \partial\Delta_n$ (i.e., $x_i=0$ if and only if $i\in I$ with some nonempty proper subset $I$ of $\{1,\ldots,n\}$) such that $F_P(x)=x$. Then,
\begin{eqnarray*}
(F_P(x))_i=\sum_{j_1,\ldots,j_{m-1}\notin I}^np_{ij_1\cdots j_{m-1}}x_{j_1}\cdots x_{j_{m-1}}=0
\end{eqnarray*}
for all $i\in I$. Hence, $p_{ij_1\cdots j_{m-1}}=0$ for all $i\in I$ and $j_1,\ldots,j_{m-1}\notin I$ which contradicts the irreducibility of $P$.

{\bf II}. $F_P$ has only finitely many fixed points in
$\mbox{rel}(\Delta_n)$.

First, by the Brouwer Fixed Point Theorem \cite{ld}, there exists at
least one fixed point in $\Delta_n$ for $F_P$. By {\bf I}, such a
fixed point is in $\mbox{rel}(\Delta_n)$. Now, suppose there are
infinitely many fixed points of $F_P$ in $\mbox{rel}(\Delta_n)$. The
compactness of $\Delta_n$ implies that there exists a sequence of
fixed points, denoted by $\{x_k\}$, such that $x_k$ converges to a
limit $\bar x\in \Delta_n$. The continuity of $F_P$ and {\bf I}
imply that $\bar x\in \mbox{rel}(\Delta_n)$. Denote by $i_d$ the
identity mapping from $\Re^n$ to itself. Then, by the assumption
that $1$ is not an eigenvalue of $\nabla F_P$ for all
$x\in\mbox{rel}(\Delta_n)$, we conclude that $i_d-F_P$ is a one to
one mapping in a small neighborhood of $\bar x$ by the Inverse
Function Theorem. However, this contradicts the fact that the
equation $(i_d-F_P)(x)$ has infinitely many solutions in the above
mentioned small neighborhood of $\bar x$ since $x_k\ra\bar x$.
Hence, $F_P$ has only finitely many fixed points in
$\mbox{rel}(\Delta_n)$.

We denote by them $\{x_1,\ldots,x_q\}$ and $U_i$ a neighborhood of
$x_i$ in $\mbox{rel}(\Delta_n)$ such that $F_P-i_d$ is a
homeomorphism between $U_i$  and a neighborhood of $0$, and injective on $\overline{ U}_i$, and the
sets $\overline{ U}_i$ are pairwise disjoint. Obviously, $F_P$ has
no fixed points in $\Delta_n\setminus(\cup_{i=1}^qU_i)$ and hence
$Ind(F_P,\mbox{rel}(\Delta_n))=\sum_{i=1}^qInd(F_P,U_i)$ by (A2) of
Theorem \ref{thm-2}.

{\bf III}. $Ind(F_P,U_i)$ is a constant for all $i\in\{1,\ldots,q\}$.

Actually, we prove a stronger result: Let $U_x$ be a small
neighborhood of $x\in \mbox{rel}(\Delta_n)$ such that $F_P-i_d$ is a
homeomorphism between $U_x$  and a neighborhood of $0$, and injective on $\overline{ U}_x$. Then
$Ind(F_P-F_P(x)+x, U_x)$ is a constant for all $x\in
\mbox{rel}(\Delta_n)$. Since $\mbox{rel}(\Delta_n)$ is relative
open, we only need to prove that $Ind(F_P-F_P(x)+x, U_x)$ is locally
constant in $\mbox{rel}(\Delta_n)$.

Now, assume that $x\in\mbox{rel}(\Delta_n)$ and $y\in U_x$. Let
$V=U_x\cap U_y$. Then $V$ is relatively open, and
\begin{eqnarray*}
Ind(F_P-F_P(y)+y,U_x)=Ind(F_P-F_P(y)+y,V)=Ind(F_P-F_P(y)+y,U_y)
\end{eqnarray*}
by the fact that $y\in V$, $F_P-i_d$ is injective on $\overline{
U}_x$ and $\overline{ U}_y$, and (A4) of Theorem \ref{thm-2}. So, it
remains to prove that
\begin{eqnarray*}
Ind(F_P-F_P(x)+x,U_x)=Ind(F_P-F_P(y)+y,U_x).
\end{eqnarray*}
To this end, define a homotopy
$H(t,\cdot):=F_P(\cdot)+t(x-F_P(x))+(1-t)(y-F_P(y))$ for $t\in[0,1]$
on $\overline{ U}_x$. Note that $H(0,\cdot)=F_P(\cdot)-F_P(y)+y$ and
$H(1,\cdot)=F_P(\cdot)-F_P(x)+x$. Now, suppose that there is
$t_0\in[0,1]$ and $x_0\in \partial U_x$ such that
$H(t_0,x_0)=F_P(x_0)+t_0(x-F_P(x))+(1-t_0)(y-F_P(y))=x_0$. By
the assumption that $x,y\in U_x$ and $F_P-i_d$ is injective on
$\overline{ U}_x$, we could conclude that $t_0\in(0,1)$ (For, if $t_0=0$, $x_0-F_P(x_0)=y-F_P(y)$, then it contradicts the fact that $F_P-i_d$ is injective on
$\overline{ U}_x$, since $x_0\in\partial U_x$ and $y\in U_x$; and, similar proof for the case when $t_0=1$). We could,
without loss of generality, shrink $(F_P-i_d)(U_x)$ to make sure
that it is a ball. So, with $x,y\in U_x$, the homeomorphism of
$F_P-i_d$ on between $U_x$ and a neighborhood of $0$, $t_0\in(0,1)$ and the Middle Value Theorem, we
could get that there exists some $x_1\in U_x$ such that
$F_p(x_1)-x_1=t_0(F_P(x)-x)+(1-t_0)(F_P(y)-y)$. So, what we get
is that $F_P(x_1)-x_1=F_P(x_0)-x_0$ for $x_1\in U_x$ and $x_0\in
\partial U_x$, which further contradicts the fact that $F_P-i_d$ is
injective on $\overline{ U}_x$. Hence, $H$ is a continuous homotopy
without fixed points on $[0,1]\times\partial U_x$. By (A3) of
Theorem \ref{thm-2}, $Ind(F_P-F_P(x)+x,U_x)=Ind(F_P-F_P(y)+y,U_x)$, which
further implies that $Ind(F_P-F_P(x)+x, U_x)$ is a constant for all $x\in
\mbox{rel}(\Delta_n)$. Now, $F_P(x_i)=x_i$ for all
$i\in\{1,\ldots,q\}$, so $Ind(F_P,U_i)$ is a constant.

{\bf IV}. $q=1$.

Choose arbitrarily $x_0\in \mbox{rel}(\Delta_n)$, and define a homotopy $H(t,\cdot):=(1-t)x_0+tF_P$ for all $t\in[0,1]$ on $\Delta_n$. Obviously, $H$ is continuous. Suppose $H(t_0,y_0)=y_0$ for some $y_0\in \partial\Delta_n$ and $t_0\in [0,1]$. By {\bf I} and $x_0\in \mbox{rel}(\Delta_n)$, $t_0\in(0,1)$. So, $y_0=(1-t_0)x_0+t_0F_P(y_0)\geq (1-t_0)x_0>0$, a contradiction to $y_0\in\partial \Delta_n$. Hence, $Ind(F_P,\mbox{rel}(\Delta_n))=Ind(x_0,\mbox{rel}(\Delta_n))$ by (A3) of Theorem \ref{thm-2}. Then, by (A1) of Theorem \ref{thm-2}, $Ind(x_0,\mbox{rel}(\Delta_n))=1$. So, $Ind(F_P,\mbox{rel}(\Delta_n))=1$. Moreover, $Ind(F_P,\mbox{rel}(\Delta_n))=\sum_{i=1}^qInd(F_P,U_i)$ by the last statement of {\bf II}.
This, together with {\bf III} which says that all $Ind(F_P,U_i)$'s are equal, implies that $q=1$.

Combining {\bf I}, {\bf II}, {\bf III} and {\bf IV}, the proof is complete. \ep

The discussion of Theorem \ref{thm-3} is motivated by \cite{k,lny}. The recent paper \cite{ln} gives a uniqueness result under some mild conditions, while they are different from the hypothesis in Theorem \ref{thm-3}.

\subsection{Third order transition probability tensors}

\hspace{4mm} We consider in this subsection specially on third order
transition probability tensors. So, it is worth describing the
considered problem more explicitly as: given a third order $n$
dimensional transition probability tensor $P$, find an $x\in\Re^n$
such that
\begin{eqnarray}\label{seq-1}
Px^2=x\;\mbox{and}\;\;x\in\Delta_n.
\end{eqnarray}

By using a result for stochastic matrices in matrix analysis \cite[Theorem 12.9]{e}, we first propose a sufficient condition to guarantee the two conditions in Theorem \ref{thm-3} for third order transition probability tensors. In the subsequent analysis, $e$ is reserved to denote the vector of all ones with appropriate size.

\begin{Theorem}\label{thm-4}
Suppose that $P$ is a third order $n$ dimensional positive transition probability tensor, and $\min_{i,j,k} p_{ijk}\geq\delta>\frac{1}{2n}$. Then, $1$ is not an eigenvalue of $\nabla F_P$ on $\Delta_n$.
\end{Theorem}

\noindent {\bf Proof.} For the convenience of the subsequent analysis, define $A_i$ as an $n\times n$ matrix with its $(j,k)$-th element being $p_{ijk}$ for all $i,j,k\in\{1,\ldots,n\}$.
So,
\begin{eqnarray}\label{eq-2}
\nabla F_P(x)=\left(\begin{array}{c}x^T(A_1+A_1^T)\\ \vdots\\x^T(A_n+A_n^T)\end{array}\right).
\end{eqnarray}

{\it {\bf Fact I} $e^T\nabla F_P(x)=2e^T$.
\begin{itemize}
\item []In fact, denote by $m_{jk}$ the $(j,k)$-th element of matrix $\nabla F_P(x)$, then
\begin{eqnarray*}
\sum_{j=1}^nm_{jk}=\sum_{j=1}^n(\sum_{i=1}^n  p_{jik}x_i+  p_{jki}x_i)=\sum_{i=1}^n(\sum_{j=1}^n  p_{jik}x_i+  p_{jki}x_i)=2\sum_{i=1}^nx_i=2.
\end{eqnarray*}
The result follows.
\end{itemize}}

By the assumption,
\begin{eqnarray}\label{eq-4}
\nabla F_P(x)\geq n\delta(\frac{2}{n}ee^T)=2\delta ee^T.
\end{eqnarray}

By {\bf Fact I}, $\frac{1}{2}\nabla F_P(x)$ is a column stochastic matrix. This, together with \reff{eq-4}, implies that $\nabla F_P(x)$ could be partitioned as
\begin{eqnarray}\label{eq-5}
\nabla F_P(x)=(1-n\delta)(2S)+n\delta(\frac{2}{n}ee^T)
\end{eqnarray}
for some column stochastic matrix $S$, here we used the fact that $x\in\Delta_n$. Actually, $\nabla F_P(x)-2\delta ee^T$ is a nonnegative matrix, and $e^T(\nabla F_P(x)-2\delta ee^T)=2e^T-2n\delta e^T=2(1-n\delta)e^T$.

{\it {\bf Fact II} Suppose that the eigenvalues of the column stochastic matrix $W$ are $\{1,\lambda_2,\ldots,\lambda_n\}$ in decreasing order of magnitude. Then, for any $\beta\in[0,1]$, the eigenvalues of the new column stochastic matrix $\beta W+(1-\beta)\frac{1}{n}ee^T$ are $\{1,\beta\lambda_2,\ldots,\beta\lambda_n\}$.

The proof is similar to that for \cite[Theorem 12.9]{e} in view of
the Perron-Frobenius Theorem \cite[Theorem 1.5]{s}.}

So, by \reff{eq-5}, {\bf Fact II} and the fact that
\begin{eqnarray*}
n\delta>\frac 12,
\end{eqnarray*}
the eigenvalues of $\nabla F_P(x)$ are $2$, $(1-n\delta)\lambda_2$, $\ldots$, $(1-n\delta)\lambda_n$ for some $\lambda_i$ with $|\lambda_i|\leq 2$.
Hence, $1$ is not an eigenvalue of $\nabla F_P(x)$. The proof is complete. \ep

\begin{Corollary}\label{cor-1}
Suppose that $P$ is a third order $n$ dimensional positive transition probability tensor, and $\sum_{j=1}^n(p_{ijk}x_j+p_{ikj}x_j)\geq\delta>\frac{1}{n}$ for all $i,j,k\in\{1,\ldots,n\}$ at $x\in\Delta_n$. Then, $1$ is not an eigenvalue of $\nabla F_P(x)$.
\end{Corollary}

\noindent {\bf Proof.} The results follows from the proof of Theorem \ref{thm-4} immediately. \ep

We note that the hypothesis in Theorem \ref{thm-4} is a sufficient
condition to guarantee that $1$ is not an eigenvalue of $F_P$ on
$\Delta_n$, hence the uniqueness of the fixed point of $F_P$ on
$\Delta_n$ by Theorem \ref{thm-3}. There may be space for refining
it. Actually, the hypothesis in Theorem \ref{thm-4} implies the
assumption in \cite[Theorem 2.3]{ln}, hence the uniqueness of the
fixed point of $F_P$ on $\Delta_n$. While, it is unknown whether the
hypothesis in \cite[Theorem 2.3]{ln} implies that $1$ is not an
eigenvalue of $F_P$ on $\Delta_n$ or not. Nevertheless, both
assumptions may be refined to guarantee the uniqueness of the fixed
point of $F_P$ on $\Delta_n$, as we could prove the following
result.

\begin{Proposition}\label{thm-1}
Let $P$ be a $2\times 2\times 2$ irreducible transition probability
tensor. Then the system \reff{seq-1} has at most one solution.
\end{Proposition}

\noindent {\bf Proof.} Let $P_{111}=\alpha\in[0,1]$. Then
$P_{211}=1-\alpha$ by the definition \reff{int-1}. Similarly, let
$P_{112}=\beta\in[0,1]$, $P_{121}=\gamma\in[0,1]$, and
$P_{122}=\tau\in[0,1]$. Then $P_{212}=1-\beta$, $P_{221}=1-\gamma$,
and $P_{222}=1-\tau$. Hence, with $s:=x_1\in[0,1]$ and $t:=x_2=1-s$,
\reff{seq-1} becomes
\begin{eqnarray}\label{seq-2}
\left\{\begin{array}{l}\alpha s^2+(\beta+\gamma)s(1-s)+\tau (1-s)^2=s,\\(1-\alpha)s^2+(2-\beta-\gamma)s(1-s)+(1-\tau)(1-s)^2=1-s.\end{array}\right.
\end{eqnarray}
System \reff{seq-2} reduces to equation
\begin{eqnarray}\label{seq-3}
(\alpha+\tau-\beta-\gamma)s^2+(\beta+\gamma-1-2\tau)s+\tau=0.
\end{eqnarray}
While, equation \reff{seq-3} has two different solutions in interval $(0,1)$ (only positive solutions, since $P$ is irreducible) is equivalent to
\begin{eqnarray}\label{seq-4}
\tau>0,\;\alpha+\tau-\beta-\gamma>0,\;\frac{1+2\tau-\beta-\gamma}{2(\alpha+\tau-\beta-\gamma)}>0,\;\mbox{and}\nonumber\\
\frac{1+2\tau-\beta-\gamma+\sqrt{(\beta+\gamma-1-2\tau)^2-4\tau(\alpha+\tau-\beta-\gamma)}}{2(\alpha+\tau-\beta-\gamma)}<1.
\end{eqnarray}
Now, equation \reff{seq-4} is equivalent to
\begin{eqnarray}\label{seq-5}
\sqrt{(\beta+\gamma-1)^2+4\tau(1-\alpha)}<2(\alpha-1)+1-\beta-\gamma.
\end{eqnarray}
Note that the right hand side of \reff{seq-5} is not greater than
$1-\beta-\gamma$, since $\alpha\in[0,1]$; while, the left hand side is
not smaller than $|1-\beta-\gamma|$, since $\tau\in[0,1]$ and
$\alpha\in[0,1]$. So, the strict inequality in \reff{seq-5} does not
hold. Hence, system \reff{seq-1} could have at most one solution.\ep

\section{Linear convergence of the power method}

\hspace{4mm} Now, we have in Theorems \ref{thm-3} and \ref{thm-4}
proved the uniqueness of the fixed point of $F_P$ on $\Delta_n$
under suitable conditions. In this section, we discuss the numerical
method for finding it and establish convergence of the method.

Here is the power method to compute a solution of system \reff{seq-1}.
\begin{Algorithm}\label{algo} (A Power Algorithm)
\begin{description}
\item [Step 0] Choose an initial guess $x^{(0)}\in\Delta_n$, let $k:=0$.

\item [Step 1] If $x^{(k)}=P(x^{(k)})^2$, stop.

\item [Step 2] Set $x^{(k+1)}:=P(x^{(k)})^2$, and $k:=k+1$. Go to Step 1.
\end{description}
\end{Algorithm}

\begin{Lemma}\label{lem-3}
Suppose that $P$ is a third order $n$ dimensional positive transition probability tensor, and $\min_{i,j,k} p_{ijk}\geq\delta>\frac{1}{2n}$. Suppose that $x^*$ is the unique fixed point of $F_P$ in $\Delta_n$, and sequence $\{x^{(k)}\}$ is generated by Algorithm \ref{algo}. Then,
\begin{eqnarray}\label{cont}
\|x^{(k+1)}-x^*\|_{1}\leq 2(1-n\delta)\|x^{(k)}-x^*\|_{1}<\|x^{(k)}-x^*\|_{1}
\end{eqnarray}
for all $k\in\{0,1,\ldots\}$. Here $\|\cdot\|_1$ means $1$-norm for vectors in $\Re^n$.
\end{Lemma}

\noindent {\bf Proof.} {\it {\bf Fact I} For the sequence $\{x^{(k)}\}$ generated by Algorithm \ref{algo}:
\begin{itemize}
\item [(a)] $x^{(k)}\in\Delta_n$ for all $k=0,1,\ldots$,
\item [(b)] $x^{(k)}_i\geq\delta$ for every $i\in\{1,\ldots,n\}$ and all $k=1,2,\ldots$.
\end{itemize}
}
Define $A_i$'s as those in the proof of Theorem \ref{thm-4}. Hence,
\begin{eqnarray}\label{eq-7}
x^{(k+1)}-x^*&=&P(x^{(k)})^2-x^*\nonumber\\
&=&P(x^{(k)})^2-P(x^*)^2\nonumber\\
&=&\left(\begin{array}{c}(x^{(k)})^TA_1x^{(k)}\\ \vdots\\(x^{(k)})^TA_nx^{(k)}\end{array}\right)-\left(\begin{array}{c}(x^{*})^TA_1x^{*}\\ \vdots\\(x^{*})^TA_nx^{*}\end{array}\right)\nonumber\\
&=&\left(\begin{array}{c}\left[(x^{(k)})^TA_1+(x^*)^TA_1^T\right](x^{(k)}-x^*)\\\vdots\\ \left[(x^{(k)})^TA_n+(x^*)^TA_n^T\right](x^{(k)}-x^*)\end{array}\right)\nonumber\\
&=&\left(\begin{array}{c}(x^{(k)})^TA_1+(x^*)^TA_1^T\\ \vdots\\(x^{(k)})^TA_n+(x^*)^TA_n^T\end{array}\right)\left(x^{(k)}-x^*\right).
\end{eqnarray}
Denote by
\begin{eqnarray*}
K:=\left(\begin{array}{c}(x^{(k)})^TA_1+(x^*)^TA_1^T\\ \vdots\\(x^{(k)})^TA_n+(x^*)^TA_n^T\end{array}\right).
\end{eqnarray*}

{\it {\bf Fact II} $e^TK=2e^T$.
\begin{itemize}
\item []In fact,  if the $(j,k)$-th element of matrix $K$ is denoted by $K_{jk}$, then
\begin{eqnarray*}
\small \sum_{j=1}^nK_{jk}=\sum_{j=1}^n(\sum_{i=1}^n  p_{jik}x^{(k)}_i+  p_{jki}x^*_i)=\sum_{i=1}^n(\sum_{j=1}^n  p_{jik}x^{(k)}_i+  p_{jki}x^*_i)=\sum_{i=1}^nx^{(k)}_i+\sum_{i=1}^nx^*_i=2.
\end{eqnarray*}
The result follows.
\end{itemize}}

By similar proof to that for Theorem \ref{thm-4}, we could get that
\begin{eqnarray*}
K=(1-n\delta)(2S)+n\delta(\frac{2}{n}ee^T)
\end{eqnarray*}
for some column stochastic matrix $S$.
Hence, with \reff{eq-7}, we have
\begin{eqnarray*}
\|x^{(k+1)}-x^*\|_{1}&=&\|K(x^{(k)}-x^*)\|_1\\
&=&\|(1-n\delta)(2S)(x^{(k)}-x^*)+n\delta(\frac{2}{n}ee^T)(x^{(k)}-x^*)\|_1\\
&=&\|(1-n\delta)(2S)(x^{(k)}-x^*)\|_1\\
&\leq&2(1-n\delta)\|S\|_{1}\|x^{(k)}-x^*\|_{1}\\
&=&2(1-n\delta)\|x^{(k)}-x^*\|_{1}\\
&<&\|x^{(k)}-x^*\|_{1}
\end{eqnarray*}
since $e^Tx^{(k)}=1=e^Tx^*$, $S$ is column stochastic and $\delta>\frac{1}{2n}$. The proof is complete.\ep

As a direct consequence of Lemma \ref{lem-3}, the following result can be established easily.
\begin{Theorem}\label{thm-7}
Suppose that $P$ is a third order $n$ dimensional positive
transition probability tensor, and $\min_{i,j,k}
p_{ijk}\geq\delta>\frac{1}{2n}$. For any initial
$x^{(0)}\in\Delta_n$, Algorithm \ref{algo} either generates a set of
finitely many points $\{x^{(0)},x^{(1)},\ldots,x^{(N)}\}$ such that
$x^{(N)}$ is the unique fixed point of $F_P$ on $\Delta_n$, or
generates an infinite sequences $\{x^{(k)}\}_{k=0}^{\infty}$ such
that it globally linearly converges to the unique fixed point of
$F_P$ on $\Delta_n$.
\end{Theorem}

\section{Convergence of the second order Markov process}
\hspace{4mm} In this section, we discuss the convergence of the second order Markov chain \reff{int-2} under the same condition as that in Theorem \ref{thm-4}. We say that the second order Markov chain is convergent if sequence $\{x^{(k)}\}$ generated by \reff{int-2} converges.

The iteration \reff{int-2} represents a map from $\Re^{2n}$ to $\Re^n$, it is hard to analyze. Then, in order to use fixed point theory, we construct an auxiliary map from a space to itself first.
With the second order Markov process \reff{int-2}, we define a nonlinear map $g:\Re^{2n}\ra\Re^{2n}$ as:
\begin{eqnarray}\label{mak-2}
g(z):=\left(\begin{array}{c}Pxy\\x\end{array}\right)
\end{eqnarray}
for all $z:=(x^T,y^T)^T\in\Re^{2n}$.

So, the second order Markov process \reff{int-2} could be rewritten as
\begin{eqnarray}\label{mak-3}
\left(\begin{array}{c}x^{(k+1)}\\x^{(k)}\end{array}\right)=g\left(((x^{(k)})^T,(x^{(k-1)})^T)^T\right)
\end{eqnarray}
for all $k=1,2,\ldots$. We could pair sequence $\{x^{(k)}\}$ successively into another sequence $\{z^{(k)}\}$ with
\begin{eqnarray}\label{mak-4}
z^{(k)}:=((x^{(k)})^T,(x^{(k-1)})^T)
\end{eqnarray}
for all $k=1,2,\ldots$. Then, the second order Markov process \reff{int-2} (equivalently \reff{mak-3}) could be further rewritten in a more compact form as
\begin{eqnarray}\label{mak-5}
z^{(k+1)}=g(z^{(k)})
\end{eqnarray}
for all $k=1,2,\ldots$.

The following lemma is straightforward.
\begin{Lemma}\label{lem-2}
For any third order $n$ dimensional transition probability tensor $P$, denote the nonlinear map associated to the second order Markov chain it induced as $g$. Then, the second order Markov process it induced as \reff{int-2} converges if and only if the sequence $\{z^{(k)}\}$ produced by $g$ as \reff{mak-5} converges.
\end{Lemma}

\begin{Theorem}\label{thm-8}
Suppose that $P$ is a third order $n$ dimensional positive transition probability tensor, and $\min_{i,j,k} p_{ijk}\geq\delta>\frac{1}{2n}$. Then, the nonlinear map $g$ has a unique fixed point $z^*\in\Delta_{n}\times\Delta_n$ and the sequence $\{z^{(k)}\}$ generated by \reff{mak-5} with initial $\Delta_{n}\times\Delta_n\ni z^{(1)}\geq \delta e>\frac{1}{2n}e$ converges globally $R$-linearly to $z^*$ as follows:
\begin{eqnarray}\label{mak-c}
\|z^{(k+1)}-z^*\|_{1}\leq (2-2n\delta)^{\lceil\frac{k+2}{2}\rceil}+(2-2n\delta)^{\lceil\frac{k+1}{2}\rceil}
\end{eqnarray}
for all $k\in\{1,2,\ldots\}$.
\end{Theorem}

\noindent {\bf Proof.} The proof consists of two parts as follows.

{\bf Part I} Obviously,
$g:\Delta_{n}\times\Delta_n\ra\Delta_{n}\times\Delta_n$ has at least
one fixed point on $\Delta_{n}\times\Delta_n$ by Theorems
\ref{thm-3} and \ref{thm-4}, since $z^*:=((x^*)^T,(x^*)^T)^T$ forms a
fixed point of $g$ with the fixed point $x^*\in\Delta_n$ of $F_P$.
Suppose $z:=(x^T,y^T)^T$ is a fixed point of $g$ on
$\Delta_{n}\times\Delta_n$. Then by \reff{mak-2},
\begin{eqnarray*}
x=Pxy,\;\mbox{and}\;\;y=x.
\end{eqnarray*}
Hence, $x\in\Delta_n$ and $x=Px^2$. So, $x$ and hence $y$ is unique by Theorems \ref{thm-3} and \ref{thm-4}.

{\bf Part II} Denote by $z^*$ the unique fixed point of $g$ on $\Delta_{n}\times\Delta_n$. We know from {\bf Part I} that $z^*=((x^*)^T,(x^*)^T)^T$ with $x^*$ the unique fixed point of $F_P$ on $\Delta_n$. Define $A_i$ as an $n\times n$ matrix with its $(j,k)$-th element being $p_{ijk}$ for all $i,j,k\in\{1,\ldots,n\}$. Denote by $((x^{(k)})^T,(y^{(k)})^T)^T:=z^{(k)}$, we have
\begin{eqnarray}\label{mak-6}
z^{(k+1)}-z^*&=&\left(\begin{array}{c}Px^{(k)}y^{(k)}-P(x^*)^2\\y^{(k+1)}-x^*\end{array}\right)\nonumber\\
&=&\left(\begin{array}{c}Px^{(k)}y^{(k)}-P(x^*)^2\\x^{(k)}-x^*\end{array}\right)\nonumber\\
&=&\left(\begin{array}{c}(y^{(k)})^TA_1^T\\ \vdots\\(y^{(k)})^TA_n^T\\I\end{array}\right)(x^{(k)}-x^*)+\left(\begin{array}{c}(x^*)^TA_1\\ \vdots\\(x^*)^TA_n\\0\end{array}\right)(y^{(k)}-x^*),
\end{eqnarray}
where we used the fact that $y^{(k+1)}=x^{(k)}$.
While
\begin{eqnarray}\label{mak-7}
z^{(k+1)}-z^*=\left(\begin{array}{c}x^{(k+1)}-x^*\\y^{(k+1)}-y^*\end{array}\right)=\left(\begin{array}{c}x^{(k+1)}-x^*\\x^{(k)}-x^*\end{array}\right).
\end{eqnarray}

By the fact that $y^{(k)}=x^{(k-1)}$, \reff{mak-6} and \reff{mak-7}, we get that
\begin{eqnarray}\label{mak-8}
x^{(k+1)}-x^*=\left(\begin{array}{c}(x^{(k-1)})^TA_1^T\\ \vdots\\(x^{(k-1)})^TA_n^T\end{array}\right)(x^{(k)}-x^*)+\left(\begin{array}{c}(x^*)^TA_1\\ \vdots\\(x^*)^TA_n\end{array}\right)(x^{(k-1)}-x^*).
\end{eqnarray}
We could prove recursively that $x^{(k)}\in\Delta_n$ and
$y^{(k)}\in\Delta_n$ for all $k=1,2,\ldots$. So, by a proof similar
to the proof of Lemma \ref{lem-3}, we could show that
\begin{eqnarray}\label{mak-9}
\|x^{(k+1)}-x^*\|_1\leq(1-n\delta)\left(\|x^{(k)}-x^*\|_1+\|x^{(k-1)}-x^*\|_1\right)
\end{eqnarray}
for all $k=1,2,\ldots$.
Now, as $x^{(1)}\geq\delta e>\frac{1}{2n}e$ and $y^{(1)}\geq\delta e>\frac{1}{2n}e$,
\begin{eqnarray*}
\|x^{(1)}-x^*\|_1\leq2-2n\delta<1,\;\mbox{and}\;\|x^{(0)}-x^*\|_1\leq 2-2n\delta<1
\end{eqnarray*}
since $x^*\geq \delta e>\frac{1}{2n}e$ and $x^{(0)}=y^{(1)}$. Hence,
\begin{eqnarray*}
\|x^{(2)}-x^*\|_1&\leq&(1-n\delta)\left(\|x^{(1)}-x^*\|_1+\|x^{(0)}-x^*\|_1\right)\leq(2-2n\delta)^2,\\
\|x^{(3)}-x^*\|_1&\leq&(1-n\delta)\left(\|x^{(2)}-x^*\|_1+\|x^{(1)}-x^*\|_1\right)\\
&\leq& (1-n\delta)\left[(2-2n\delta)^2+(2-2n\delta)\right]\leq(2-2n\delta)^2.
\end{eqnarray*}
So, inductively, we could prove that
\begin{eqnarray*}
\begin{array}{c}
\|x^{(k+1)}-x^*\|_{1}\leq (2-2n\delta)^{\lceil\frac{k+2}{2}\rceil},\\
\|y^{(k+1)}-x^*\|_{1}=\|x^{(k)}-x^*\|_{1}\leq(2-2n\delta)^{\lceil\frac{k+1}{2}\rceil},\end{array}
\end{eqnarray*}
for all $k=1,2,\ldots$.
The proof is complete.  \ep

\section{Numerical examples}
\hspace{4mm} In this section, we present some numerical examples to show the feasibility of the results in Section 4. The first numerical example is taken from \cite{ln}.

\begin{Example}\label{exm-1}
{\em This example comes from DNA sequence data in Tables 6 and 10
of \cite{rt}. There are two third order three dimensional transition probability tensors. By using
the MatLab multi-dimensional array notation, the transition transition probability tensors are
given by
\begin{eqnarray*}
(i)\begin{array}{c}P(:,:,1)=\left(\begin{array}{ccc}0.6000& 0.4083& 0.4935\\0.2000& 0.2568& 0.2426\\0.2000 &0.3349& 0.2639\end{array}\right),\;
P(:,:,2)=\left(\begin{array}{ccc}0.5217&0.3300&0.4152\\0.2232&0.2800&0.2658\\0.2551&0.3900&0.3190\end{array}\right),\\
P(:,:,3)=\left(\begin{array}{ccc}0.5565& 0.3648& 0.4500\\0.2174& 0.2742& 0.2600\\0.2261& 0.3610& 0.2900\end{array}\right),\end{array}
\end{eqnarray*}
and
\begin{eqnarray*}
(ii)\begin{array}{c}P(:,:,1)=\left(\begin{array}{ccc}0.5200& 0.2986& 0.4462\\0.2700& 0.3930& 0.3192\\0.2100& 0.3084 &0.2346\end{array}\right),\;
P(:,:,2)=\left(\begin{array}{ccc}0.6514 &0.4300 &0.5776\\0.1970& 0.3200& 0.2462\\0.1516& 0.2500 &0.1762\end{array}\right),\\
P(:,:,3)=\left(\begin{array}{ccc}0.5638 &0.3424 &0.4900\\0.2408 &0.3638& 0.2900\\0.1954 &0.2938& 0.2200\end{array}\right),\end{array}
\end{eqnarray*}
respectively.

Note that tensor $(i)$ satisfies the assumption in Theorem
\ref{thm-4}. We compute the iteration sequence through both the second order power method
(i.e., Algorithm \ref{algo}) and the Markov process iteration
\reff{int-2}. For every case, we choose randomly the initial guess
$x^{(0)}\in\Delta_n$ for the power method and use $x^{(0)}$ and
$x^{(1)}:=P(x^{(0)})^2$ for the Markov iteration. The corresponding
algorithm is terminated whenever $\|x^{(k)}-x^{(k-1)}\|_1<1.0\times
10^{-6}$. We simulate ten times for every case and record the
average number of iterations {\bf It} and {\bf Itm} for the power
method and the Markov iteration, respectively. For (i), ${\bf
It}=9.1$ and ${\bf Itm}=12$; and for (ii), ${\bf It}=5.8$ and ${\bf
Itm}=11$. The computed stationary probability vector in every
simulation coincides with that in \cite[Example 1]{ln}. In order to illustrate the convergent rates in the above sections. We pictured the values of $\frac{\|x^{(k+1)}-x^*\|_1}{\|x^{(k)}-x^*\|_1}$ corresponding to the Power method ({\bf Power}) of one random test in Figure 1, and the values of $\|x^{(k)}-x^*\|_1$ corresponding to the Markov process iteration ({\bf Markov}) of one random test in Figure 2. The theoretical bounds ({\bf Theoretical}) for both methods (Lemma \ref{lem-3} and Theorem \ref{thm-8}) are pictured in the corresponding figures as well. Figure 1 demonstrates the global linear convergence, the last point being zero is due to that we choose $x^*$ to be the last iteration. Figure 2 demonstrates the global $R$-linear convergence.}
\end{Example}
\begin{Example}\label{exm-2}
{\em Third order probability tensors with dimensions $n=100$ are generated in this example. The details are: a positive tensor is generated randomly with its every entry in $(0,1)$, and then scale the resulting tensor to be a probability tensor. Add every entry of the probability tensor with $\frac{\delta}{1-n\delta}$, where $\delta:=\frac{13}{20n}$. Finally, scale the resulting tensor to be a probability tensor. We see that this tensor satisfies the assumption in Theorem \ref{thm-4} with the above $\delta$. Now, the Markov process iteration
\reff{int-2} is used to compute the stationary probability vector. We simulated ten times and pictured the values of $\|x^{(k)}-x^*\|_1$ for every simulation in Figure 3. The theoretical bound in Theorem \ref{thm-8} is also given. It is easy to see that the iteration curves are dominated by the theoretical curve very well.}
\end{Example}

\begin{figure}[htbp]
\centering
\includegraphics[width=4.6in]{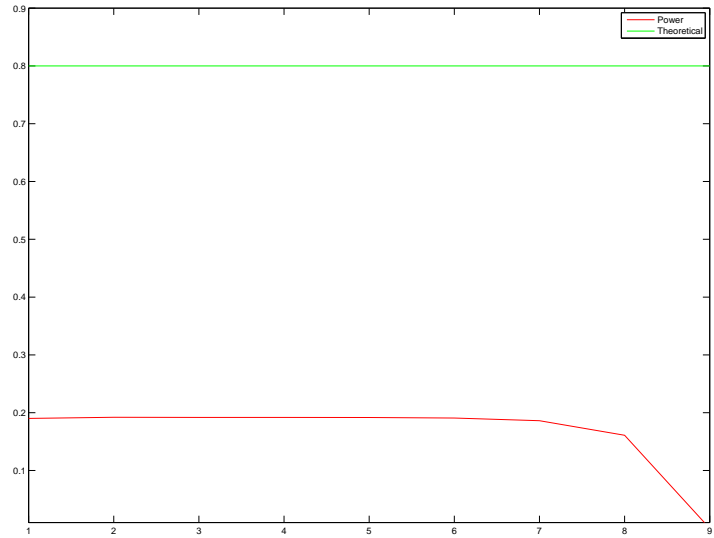}
\caption{Iteration of Example \ref{exm-1} (i): power method.}
\end{figure}

\begin{figure}[htbp]
\centering
\includegraphics[width=4.6in]{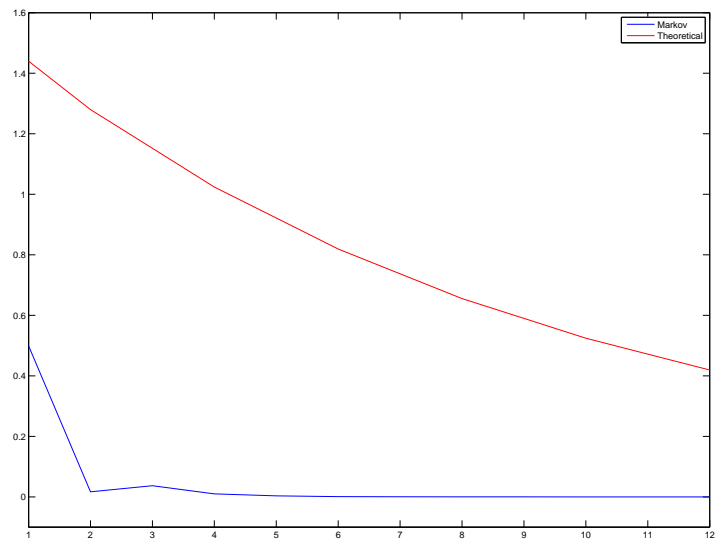}
\caption{Iteration of Example \ref{exm-1} (i): the Markov process.}
\end{figure}

\begin{figure}[htbp]
\centering
\includegraphics[width=4.6in]{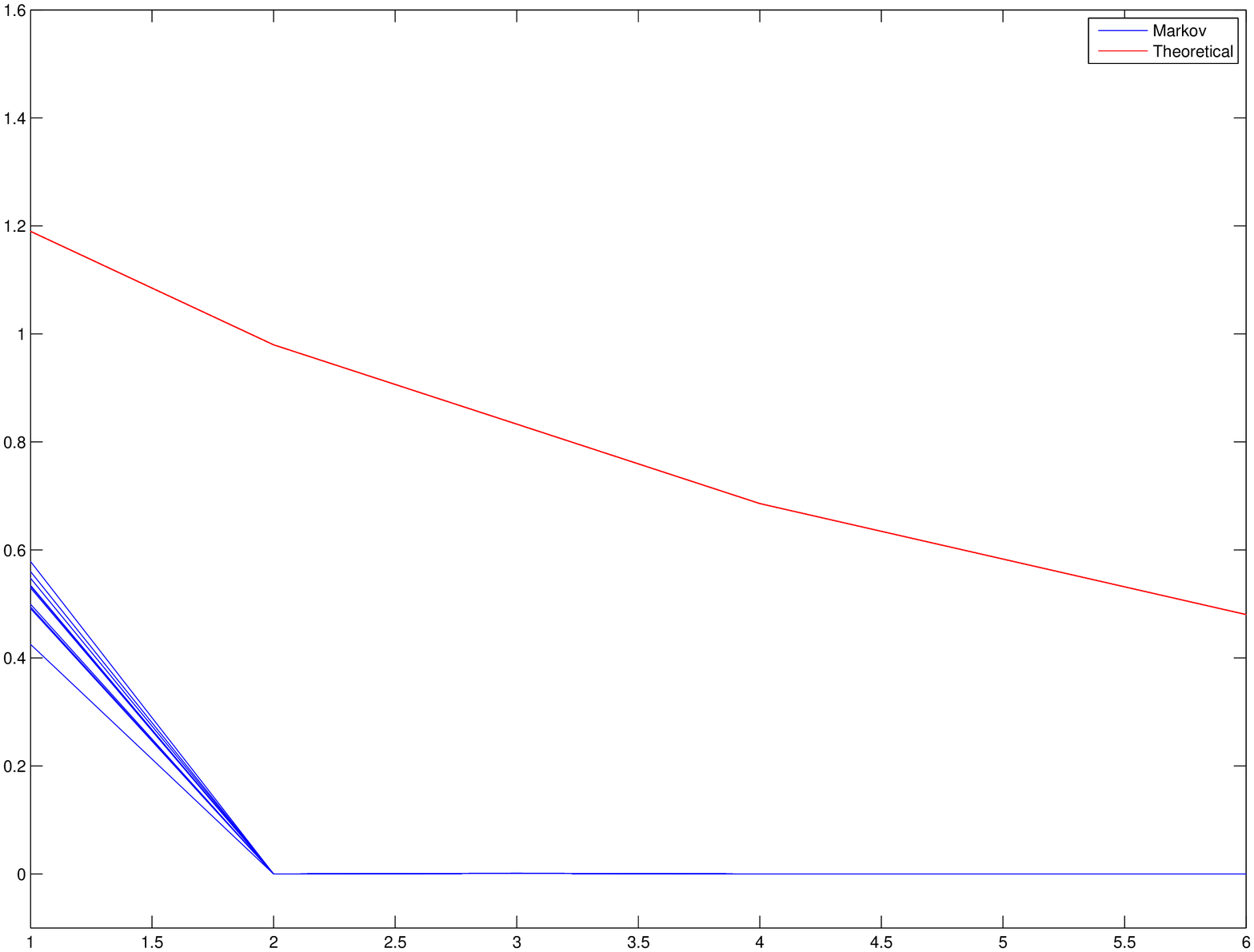}
\caption{Iteration of randomly generated tensors with dimension $100$.}
\end{figure}

{\bf Acknowledgement.} We are grateful to Prof. Kungching Chang for
his valuable comments, Yisheng Song for suggesting us the concept of
fixed point index, and Prof. Yinyu Ye for inspiring us the nonlinear
map in Section 4.


\end{document}